\documentclass{amsart}
\newtheorem{thm}{Theorem}[section]

\newtheorem{lem}[thm]{Lemma}
\newtheorem{prop}[thm]{Proposition}
\theoremstyle{definition}
\newtheorem{defn}[thm]{Definition}
\theoremstyle{remark}

\numberwithin{equation}{section}

\newcommand{\CC}{\mathbb C}
\newcommand{\C}{\mathbb C}

\def\ll{{\mathcal L}}

\begin{document}

\title[Naturally graded 2-filiform Leibniz algebras.]
{Naturally graded 2-filiform Leibniz algebras.}%
\author{L.M. Camacho, J. R. G\'{o}mez}
\address{[L.M. Camacho -- J.R. G\'{o}mez] Dpto. Matem\'{a}tica Aplicada I. Universidad de Sevilla. Avda. Reina Mercedes, s/n. 41012 Sevilla. (Spain)}
\email{lcamacho@us.es --- jrgomez@us.es}

\author{A.J. Gonz\'{a}lez}
\address{[A.J. Gonz\'{a}lez] Dpto. de Matem\'{a}ticas. Universidad de Extremadura. Avda. de la Universidad, s/n. Badajoz. (Spain)}
\email{agonzale@unex.es}

\author{B.A. Omirov}
\address{[B.A. Omirov] Institute of Mathematics of Academy of Uzbekistan. 29, F.Khodjaev str. 700125, Tashkent, Uzbekistan}
\email{omirovb@mail.ru}



\thanks{The last author was supported by grant of INTAS Ref. Nr. 0483-3035 and he is very gratefull
to hospitality the Universidad de Sevilla}%

\begin{abstract}
 The Leibniz algebras appear as a generalization of the Lie
algebras \cite{loday}. The classification of naturally graded $p$-filiform Lie
algebras is known \cite{C-G-JM}, \cite{J.Lie.Theory}, \cite{AJM}, \cite{Ve}. In
this work we deal with the classification of $2$-filiform Leibniz
algebras. The study of $p$-filiform Leibniz non Lie
algebras is solved for $p=0$ (trivial) and $p=1$
\cite{Omirov1}. In this work we get the classification of
naturally graded non Lie $2$-filiform Leibniz algebras.
\end{abstract}
\maketitle
\section{Introduction}

In this work we study the naturally graded 2-filiform Leibniz
algebras. Since the filiform (1-filiform) Lie algebras have the
maximal nilindex, Vergne studied them and obtained the
classification of naturally graded \cite{Ve}. Many authors have
studied the complete classification for low dimensions. The
lists up to dimension 8 can be found in \cite{G-H} and the
classification filiform up to dimension 11 in \cite{b6} The
notion of $p$-filiform Lie (resp. Leibniz) algebras can be
considered as a generalization of filiform Lie algebras.

The knowledge of naturally graded algebras of a certain family
offers significant information about the complete family. The
classification of $2$-filiform Lie algebras and $p$-filiform
has been obtained \cite{AJM}, \cite{J.Lie.Theory}.

In the case of Leibniz algebras only the classification of
$0$-filiform and $1$-filiform algebras is known \cite{Omirov1},
\cite{Omirov2}. In the present paper we get the classification of
naturally graded $2$-filiform Leibniz algebras.

Leibniz algebras are defined by the Leibniz identity:
$$[x,[y,z]]=[[x,y],z]-[[x,z],y]$$

We have used the software {\it Mathematica} to study particular cases
in concrete finite dimensions and later, by induction,
the obtained results are generalized for arbitrary
finite dimension.

Let $\ll$ be a Leibniz algebra, we define the following sequence:
$$\ll^1=\ll,\quad \ll^{n+1}=[\ll^n,\ll]$$

An algebra $\ll$ is {\it nilpotent} if $\ll^n=0$ for some $n \in {\bf N}$.

For any element $x$ of $\ll$ we define $R_x$ the operator of right multiplication as
$$\begin{array}{llll}
R_x:&z&\rightarrow & [z,x], \ z \in \ll
\end{array}$$

Let us $x\in \ll \backslash [\ll,\ll]$ and for the nilpotent operator
$R_x$ of right multiplication, define the decreasing sequence
$C(x)=(n_1,n_2,...,n_k)$ that consists of the dimensions of the
Jordan blocks of the $R_x$. Endow the set of these sequences with the
lexicographic order.

The sequence $C(\ll)=max_{x\in \ll \backslash [\ll,\ll]} C(x)$ is defined to be the
{\it characteristic sequence} of the algebra $\ll$.

\begin{defn}

A Leibniz algebra $\ll$ is called $p$-filiform if
$C(L)=(n-p,\underbrace{1,\dots,1}_{p})$, where $p\geq 0$.
\end{defn}

Note that this definition when $p>0$ agrees with the definition
of $p$-filiform Lie algebras.

From now we will use the expression ``graded algebra" instead of
``naturally graded algebra".

Let $\ll$ be a graded
$p$-filiform $n$-dimensional Leibniz algebra, then there exists
a basis $\{e_1, e_2, \dots,e_n\}$ such that $e_1\in \ll-\ll^2$ and
$C(e_1)=(n-p,\underbrace{1,\dots,1}_{p})$.

By definition of
characteristic sequence the operator $R_{e_1}$ in Jordan form has one
block $J_{n-p}$ of size $n-p$ and $p$ block $J_1$
(where $J_1=\{0\}$) of size one.

The possibilities for operator $R_{e_1}$ are the follow:
$$\small\begin{array}{llll}
\left(\begin{array}{ccccc}
J_{n-p}& 0&0 & \cdots& 0\\
0&J_1&0&\cdots&0\\
\vdots& & &\vdots&\\
0&0&0&\cdots& J_1
\end{array}\right),& \left(\begin{array}{ccccc}
J_{1}& 0&0 & \cdots& 0\\
0&J_{n-p}&0&\cdots&0\\
\vdots& & &\vdots&\\
0&0&0&\cdots& J_1
\end{array}\right),&\cdots,& \\[10mm]
\left(\begin{array}{ccccc}
J_{1}& 0&0 & \cdots& 0\\
0&J_1&0&\cdots&0\\
\vdots& & &\vdots&\\
0&0&0&\cdots& J_{n-p}
\end{array}\right)& &
\end{array}$$

It is easy to prove that when $J_{n-p}$ is placed an a different position from the first
are isomorphic cases. Thus,
we have only the following possibilities of Jordan form of the matrix $R_{e_1}$:
$$\small\begin{array}{ll}
\left(\begin{array}{ccccc}
J_{n-p}& 0&0 & \cdots& 0\\
0&J_1&0&\cdots&0\\
\vdots& & &\vdots&\\
0&0&0&\cdots& J_1
\end{array}\right),& \left(\begin{array}{ccccc}
J_{1}& 0&0 & \cdots& 0\\
0&J_{n-p}&0&\cdots&0\\
\vdots& & &\vdots&\\
0&0&0&\cdots& J_1
\end{array}\right)
\end{array}$$

\begin{defn}
A $p$-filiform Leibniz algebra $\ll$ is called first type (type I) if the operator $R_{e_1}$ has the form:
$$\left(\begin{array}{ccccc}
J_{n-p}& 0&0 & \cdots& 0\\
0&J_1&0&\cdots&0\\
\vdots& & &\vdots&\\
0&0&0&\cdots& J_1
\end{array}\right)$$
and second type (type II) in the other case.
\end{defn}

\subsection{Naturally graded filiform and $2$-filiform Lie algebras}

Naturally graded $p$-filiform Lie algebras are known for all $p>0$, \cite{J.Lie.Theory}, \cite{AJM},  \cite{Ve}.

\vglue3mm

Examples of filiform Lie algebras are $\ll_n$, $Q_n$ defined as follows:

$$
\begin{array}{l}
\begin{array}{l}
{\mathcal L}_n\quad (n\geq 3): \left\{ \ [X_0,X_i]=X_{i+1} \ \quad 1\leq i\leq n-2. \right.
\end{array}
\\[6mm]

\begin{array}{l}
{\mathcal Q}_n\quad (n\geq 6, \ n\ {\rm even}):\left\{
\begin{array}{ll}
[X_0,X_i]=X_{i+1} & \quad 1\leq i\leq n-2, \\[2mm]
[X_i,X_{n-1-i}]=(-1)^{i-1}\,X_{n-1} & \quad 1\leq i \leq
\frac{n-2}{2}.
\end{array}
\right.

\end{array}
\end{array}
$$

Provided examples of $2$-filiform Lie algebras.

\vglue3mm
$ {\mathcal L}(n,r) \quad
\left(n\geq 5, \
3\leq r\leq 2\,\left\lfloor\frac{n-1}{2}\right\rfloor -1, \ r
\  {\rm odd}
\right):
$

\vglue1mm
%
$\qquad
\left\{
\begin{array}{p{6cm}l}
$  [X_0, X_i]= X_{i+1} $ & 1\leq i\leq n-3, \\[3mm]
$  [X_i, X_{r-i}]= (-1)^{i-1}\,Y $ &  1\leq i
 \leq\frac{r-1}{2}.
\end{array}
\right.
$

\vglue3mm
${\mathcal Q}(n,r) \quad
\left(
 n\geq7,\ n \mbox{\rm \ odd}; \
        3\leq r\leq n-4, \ r\ {\rm odd}
\right)
$:

\vglue1mm
$\qquad
\left\{
\begin{array}{p{6cm}l}
$ [X_0, X_i]= X_{i+1} $ & 1\leq i\leq n-3, \\[3mm]
$ [X_i, X_{r-i}]= (-1)^{i-1}\,Y $& 1\leq i\leq\frac{r-1}{2},\\[3mm]
$ [X_i, X_{n-2-i}]= (-1)^{i-1}\, X_{n-2}$ & 1\leq i\leq\frac{n-3}{2}.
\end{array}
\right.
$

\vglue3mm
$\tau(n,n-4)$ \ ($n$ odd, $n\geq 7$):

\vglue1mm
$\qquad
\left\{
\begin{array}{p{8.6cm}l}
$ [X_0, X_i]= X_{i+1} $ &  1\leq i\leq n-3, \\{}
$ [X_i, X_{n-4-i}]= (-1)^{i-1}\,(X_{n-4}+Y)$ &
     1\leq i\leq\frac{n-5}{2},\\{}
$  [X_i, X_{n-3-i}]= (-1)^{i-1}\,\frac{(n-3-2i)}{2}\,X_{n-3}$ &
                1\leq i\leq\frac{n-5}{2}, \\{}
$  [X_i, X_{n-2-i}]= (-1)^{i}\,(i-1)\,\frac{(n-3-i)}{2}\,X_{n-2} $&
               2\leq i\leq\frac{n-3}{2}, \\{}
$  [X_i,Y]= \frac{(5-n)}{2}\,X_{n-4+i} $ &   1\leq
i\leq2.
\end{array}
\right.
$

\vglue3mm
$\tau(n,n-3)$ \ ($n$  even, $n\geq 6$):

\vglue1mm
$\qquad
\left\{
\begin{array}{p{8.6cm}l}
$  [X_0, X_i]= X_{i+1} $ &   1\leq i\leq n-3, \\{}
$  [X_i, X_{n-3-i}]= (-1)^{i-1}\,(X_{n-3}+Y) $ &
       1\leq i\leq\frac{n-4}{2},\\{}
$  [X_i, X_{n-2-i}]= (-1)^{i-1} \, \frac{(n-2-2i)}{2}\,X_{n-2} $ &
         1\leq i\leq\frac{n-4}{2}, \\{}
$  [X_1,Y]= \frac{(4-n)}{2}\,X_{n-2}. $ &
\end{array}
\right.
$

\section{Naturally graded $p$-filiform Leibniz algebra}

It is easy to see that a Leibniz algebra of type I is not a Lie algebra.

Let $\ll$ be an $n$-dimensional $p$-filiform Leibniz algebra. We
define a natural gradation of $\ll$ as follows. Take $\ll_1=\ll,$
$\ll_i=\ll^i/\ll^{i+1},\ 2\leq i\leq n-p$. It is clear that $\ll
\simeq \\_1\oplus\ll_2\oplus \cdots \oplus \ll_{n-p},$ where
$[\ll_i,\ll_j]\subseteq \ll_{i+j}$ and $\ll_{i+1}=[\ll_i,\ll_1]$
for all $i$.

Let $\ll$ be a graded $p$-filiform Leibniz algebra of the first type.
Then there exists a basis $\{e_1, e_2, \dots, e_{n-p},f_1,\dots,f_p\}$ such that
$$\begin{array}{ll}
[e_i, e_1]=e_{i+1},& 1 \leq i \leq n-p-1\\{}
[f_j, e_1]=0,&  1 \leq j \leq p.
\end{array}$$

From this multiplication we have: $$<e_1> \subseteq  \ll_1,\quad <e_2>
\subseteq \ll_2,\quad <e_3> \subseteq \ll_3,\dots, <e_{n-p}> \subseteq
\ll_{n-p}$$ but we do not know about the places of the elements
$\{f_{1}, f_{2},\dots,f_p\}$.

\

Let denote by $r_1, r_2,\dots, r_p$ the places of elements $f_{1},$ $ f_{2},$ $\dots,$ $f_p$ in natural
gradation  correspondingly, that is,
$f_{i}\in \ll_{r_i}$ with $1\leq i \leq p$. Further the law of a
Leibniz algebra of type I with the set $\{r_1, r_2,\dots,r_p\}$ will be denoted by
$\mu_{(I,r_1,\dots,r_p)}$.

\

We can suppose that $1 \leq r_1 \leq r_2 \leq \cdots \leq r_p \leq n-p$.

\

\begin{thm}
 Let $\ll$ be a graded $p$-filiform Leibniz algebra of type I. Then $r_s \leq s$
 for any $s\in \{1, 2,\dots, p\}$.
\end{thm}

{\bf Proof:}

Note $r_1=1$. In fact, if $r_1 > 1$, then the algebra $\ll$ is
one generated and by [\cite{Omirov1}, lemma $1$] it is nul-filiform Leibniz
algebra, and hence $C(\ll)=(n, 0)$, that is, we obtain contradiction with
condition $C(L)=(n-p, 1, 1,\dots, 1)$.

\

Let us prove that $r_2 \leq 2$.
Suppose otherwise, that is, $r_2 > 2$. Then $$\ll_1=<e_1, e_{n-p+1}>$$ $$\ll_2=<e_2>$$
$$\ll_{r_2}=[\ll_{r_2-1},\ll_1]=<[<e_{r_2-1}>,<e_1,f_{1}>]>=<e_{r_2},[e_{r_2-1} ,f_{1}]>$$

\

Consider the multiplication:
$$[e_{r_2-1},f_{1}]=[[e_{r_2-2} ,e_1], f_{1}]=[ e_{r_2-2}, [e_1, f_{1}]] + [[e_{r_2-2} , f_{1}], e_1]$$

\

Since the multiplication $[e_1, f_{1}]\in \ll_2=<e_2>\subseteq
Z(\ll)$, then the first item is equal to zero. It is evident that the
second item belongs to the linear span  $<e_{r_2} >$.  So,
$f_{2}\notin \ll_{r_2}$ and we obtain contradiction method, hence
$r_2\leq 2$.

\

Let us suppose that the condition of the theorem is true for any
value less than $s$. We prove that $r_s \leq s$. We shall prove
it by contradiction, that is, suppose that $r_s>s$.

\

If $r_s  > s$ we prove the following embedding:
$$[e_{r_s-r_t} ,f_{t}]\subseteq <e_{r_s}>,\qquad 1 \leq  t \leq s-1$$

\

We shall prove it by descending induction by t.

\

Let us prove it for $t=s-1$. Consider the multiplication:
$$\begin{array}{lll}
[e_{r_s-r_{s-1}},f_{s-1}]&=&[[e_{r_s-r_{s-1}-1},e_1],f_{s-1}]=[e_{r_s-r_{s-1}-1},[e_1,f_{s-1}]]+\\
&+&[[e_{r_s-r_{s-1}-1},f_{s-1}],e_1]
\end{array}$$

\

Since $r_s  > s$, we have $[e_1, f_{s-1}]\in \ll_{r_{s-1}+1}
=<e_{r_{s-1}+1} >\in Z(\ll)$, that is,
$[e_{r_s-r_{s-1}-1},[e_1,f_{s-1}]]=0$. From the multiplication on the right side on $e_1$
we have
$$[[e_{r_s-r_{s-1}-1},f_{s-1}],e_1]\subseteq <e_{r_s}>$$
hence, $[e_{r_s-r_{s-1}},f_{s-1}]\subseteq <e_{r_s}>$.

\

Let suppose that embedding $[e_{r_s-r_t},f_{t}]\subseteq <e_{r_s}>$ is true for any value
greater than $t+1$. We prove it for $t$.

\

Consider the multiplication:
$$\begin{array}{l}
[e_{r_s-r_t},f_{t}]=[[e_{r_s-r_t-1},e_1],f_{t}]=[e_{r_s-r_t-1},[e_1,f_{t}]]+\\
+[[e_{r_s-r_t-1},f_{t}],e_1]
\end{array}$$

As $[e_1, f_{t}]\in \ll_{r_t+1}$, then in case
$r_t+1=r_{t+1}$ we have $\ll_{r_t+1}=\{e_{r_{t+1}}, f_{t+1} \
\vee  \  f_{t+2}\ \vee \cdots \vee \ f_{s-1}\}$. Therefore
the multiplication $[e_{r_s-r_t-1},[e_1,f_{t}]]$  is contained
in linear span $<e_{r_s}>$  by induction.

\

If $r_t+1\neq r_{t+1}$ the following equality $[e_1,
f_{t}]=<e_{r_t+1} >$ is hold (because $\ll_{r_t+1}
=<e_{r_t+1}>$) and hence $[e_{r_s-r_t-1},[e_1,f_{t}]]=0$. Evidently, the
second item also is contained in the linear span $<e_{r_s}>$.

\

Thus, $[e_{r_s-r_t},f_{t}]\subseteq <e_{r_s}>,$ \ $1\leq t\leq s-1$
is proved.

Let us prove that  $\ll_{r_s}\subseteq <e_{r_s}>$ supposing $r_s>s$.
Consider the multiplication:
$$\ll_{r_s}=[\ll_{r_s-1},\ll_1]=[<e_{r_s-1}>,<e_1,f_{1}\vee \cdots \vee f_{s-1}>]$$

From $[e_{r_s-r_t},f_{t}]\subseteq
<e_{r_s}>$, we have that $\ll_{r_s}\subseteq <e_{r_s}>$,
that is, we obtain the contradiction which completes the proof of
theorem.

$\hfill\Box$

\subsection{Naturally graded $2$-filiform Leibniz algebras}

In this section naturally graded $2$-filiform Leibniz algebras will be classified.

The classification of the null-filiform Leibniz algebras is an easy task
and one for naturally graded $1$-filiform Leibniz algebras is
similar to the case of Lie algebras. However, when $p$, increases
the difficulties also increase exponentially in the study of Leibniz
algebras with respect to Lie algebras.

From \cite{AJM} we observe the existence of graded
$2$-filiform Lie algebras in arbitrary dimension. Let us
demonstrate examples of graded $2$-filiform Leibniz algebras of
type I which obviously are not Lie algebras.

In this work, we use the following notation:
\begin{itemize}
\item $\{e_1,e_2,...,e_{n-2},e_{n-1},e_n\}$ an adapted basis and
\item $r_1,\ r_2$ the places of elements $e_{n-1},\ e_n$.
\end{itemize}

{\bf Example 1.} Let $\ll_{n-2}^0$  be a graded nul-filiform
Leibniz algebra of dimension $n-2$ and $\ll_{n-1}^1$  be a graded
filiform non Lie Leibniz algebra of dimension $n-1$ of type I. Then
$\ll_{n-2}^0\oplus {\bf C}^2$ and  $\ll_{n-1}^1\oplus {\bf C}$
are graded n-dimensional split $2$-filiform Leibniz algebras of
type I.

\

The following lemma establishes that a graded $2$-filiform
Leibniz algebra of type I with condition $r_1=r_2=1$ is a split algebra
from the above example.

\begin{lem}\label{(I,1,1)}
Let $\ll$ be a graded $2$-filiform
Leibniz algebra of type  $\mu_{(I, 1, 1)}$. Then $\ll$ is a split algebra from
example 1.
\end{lem}

{\bf  Proof:}

Let algebra $\ll$ has form  $\mu_{(I, 1, 1)}$, then for an
adapted basis $\{e_1, e_2,\dots, e_n\}$ the multiplications on the right
side on $e_1$ are the following:

$$\left\{\begin{array}{ll}
[e_i, e_1]=e_{i+1}, &  1\leq i\leq n-3 \\[1mm]
[e_i, e_{n-1}]={\alpha_i}e_{i+1}, &  1\leq i\leq n-3  \\[1mm]
[e_{n-1}, e_{n-1}]={\alpha}_{n-1}e_2& \\[1mm]
[e_n, e_{n-1}]={\alpha}_n e_2 &\\[1mm]
[e_i, e_n]={\beta}_ie_{i+1}, & 1\leq i\leq n-3 \\[1mm]
[e_{n-1}, e_n]={\beta}_{n-1}e_2& \\[1mm]
[e_n, e_n]={\beta}_n e_2&
\end{array}\right.$$

Using Leibniz identity it is not difficult to obtain the following restrictions:
$$\left\{\begin{array}{ll}
\alpha_i=\alpha, & 1\leq i\leq n-3\\
\beta_i=\beta, & 1\leq i\leq n-3\\
\alpha_{n-1}=\alpha_n=0&\\
\beta_{n-1}=\beta_n=0
\end{array}\right.$$

Let us rewrite the multiplications of basis elements taking into account the above restrictions:
$$\left\{\begin{array}{ll}
[e_i, e_1]=e_{i+1}, &  1\leq i\leq n-3 \\[1mm]
[e_i, e_{n-1}]={\alpha}e_{i+1}, &  1\leq i\leq n-3 \\[1mm]
[e_i, e_n]={\beta}e_{i+1}, & 1\leq i\leq n-3
\end{array}\right.$$

If $\alpha\neq 0$ we take the change of basis: $e'_i=e_i,\ 1
\leq i \leq n-1,$ $e'_n = \alpha e_n-\beta e_{n-1}$, we can
suppose that the coefficient $\beta$ is equal to zero, that is, we have the
multiplications:
$$\left\{\begin{array}{ll}
[e_i, e_1]=e_{i+1}, \quad  1\leq i\leq n-3 \\[1mm]
[e_i, e_{n-1}]={\alpha}e_{i+1}, \quad  1\leq i\leq n-3
\end{array}\right.$$

If  $\alpha=0$, then taking  $e'_i=e_i,\ 1 \leq i \leq n-2$,
$e'_{n-1}=e_n,\ e'_n =e_{n-1}$, we can also suppose that
coefficient  $\beta=0$.
In this case it is easy to see that $\ll=\ll_{n-2}^0\oplus \C^2$.

\

If $\alpha\neq 0$, the change of basis
$e'_{n-1}=\displaystyle\frac{1}{\alpha}e_{n-1}$ (and $e'_i=e_i,\
i\neq n-1$) allows us to suppose $\alpha=1$ and so $\ll=\ll_{n-1}^1\oplus \CC$.

$\hfill\Box$

\

For graded non split $2$-filiform Leibniz algebra of type I with condition $r_2=2$ the following theorem is hold.

\

The next results were supported by $Mathematica$ package.

\

\begin{prop}
Let $\ll$ be an $4$-dimensional graded $2$-filiform non split Leibniz algebra of type $\mu_{(I, 1, 2)}$.Then $\ll$
is isomorphic to the following algebra:
$$\left\{\begin{array}{l}
[e_1, e_1]=e_2 \\[1mm]
[e_1, e_3]=e_4
\end{array}\right.$$
\end{prop}

{\bf Proof:}

We have that the natural gradation is:
$$<e_1,e_3>\oplus <e_2,e_4>$$
and the multiplication is:
$$\left\{\begin{array}{l}
[e_1,e_1]=e_2\\{}
[e_1,e_3]=\alpha_1 e_2+\beta_1 e_4\\{}
[e_3,e_3]=\alpha_2 e_2+\beta_2 e_4
\end{array}\right.$$
with $\beta_1\neq 0$ or  $\beta_2\neq 0$.

If we make the following change of basis $\beta_2' e'_4=\alpha_2 e_2+\beta_2 e_4$
it is possible to suppose $\alpha_2=0$ and
$$\left\{\begin{array}{l}
[e_1,e_1]=e_2\\{}
[e_1,e_3]=\alpha_1 e_2+\beta_1 e_4\\{}
[e_3,e_3]=\beta_2 e_4
\end{array}\right.$$
with $\beta_1\neq 0$ or  $\beta_2\neq 0$.

According to the characteristic sequence we have that $rank (R_{e_1+Ae_3})\leq 1$, it implies that
$\beta_2=0$ and $\beta_1\neq 0$. An elementary change of basis permits to prove this result.

$\hfill\Box$

\begin{prop}

Let $\ll$ be a $5$-dimensional naturally
graded 2-filiform Leibniz algebra of type $\mu_{(I, 1, 2)}$. Then, $\ll$ is
isomorphic to the one of the following pairwise non isomorphic
algebras:
$$\begin{array}{ll}
\mu^1:\left\{\begin{array}{ll}
[e_i, e_1]=e_{i+1}, &  1\leq i\leq 2 \\[1mm]
[e_1, e_4]=e_2+e_5, & \\[1mm]
[e_2, e_4]=e_3, &
\end{array}\right.&
\mu^2 :\left\{\begin{array}{ll}
[e_i, e_1]=e_{i+1}, &  1\leq i\leq 2 \\[1mm]
[e_1, e_4]=e_5. &
\end{array}\right.
\end{array}$$
$$\begin{array}{ll}
\mu^3:\left\{\begin{array}{ll}
[e_i, e_1]=e_{i+1}, &  1\leq i\leq 2 \\[1mm]
[e_1, e_4]=ie_2+e_5, & \\[1mm]
[e_2, e_4]=ie_3, &  i^2=-1\\[1mm]
[e_5,e_4]=e_3.
\end{array}\right.&
\mu^4: \left\{\begin{array}{ll}
[e_i, e_1]=e_{i+1}, &  1\leq i\leq 2 \\[1mm]
[e_1, e_4]=e_5.\\[1mm]
[e_5,e_4]=e_3.
\end{array}\right.
\end{array}$$
\end{prop}

{\bf Proof:}

Analogously as in above we can assume that

$$\ll_1=<e_1, e_4>,\ \ll_2=<e_2, e_5>,\ \ll_3=<e_3>, \
\ll_4=<0>$$

Put $[e_5,e_4]=\gamma e_3$. If $\gamma=0,$ then we obtain algebra
$L(\alpha,0)$, otherwise not restricted of generality we obtain
algebra $L(\alpha,1).$ Since dimension of left annihilator of the
algebra $L(\alpha,0)$ is equal to 2 ($e_4, e_5 \in L(\ll)$) and dimension of
left annihilator of the algebra $L(\alpha,1)$ is equal to 1
($e_4\in L(\ll)$) there are not isomorphic.

From the above argumentation we have following
algebras
$$\begin{array}{ll}
L(\alpha,1): \left\{\begin{array}{l} [e_{i},e_1]=e_{i+1},\quad 1\leq i\leq 2\\{}
[e_1,e_4]=\alpha e_2+f_2\\{}
 [e_2,e_4]=\alpha e_3  \\{}
 [f_2,e_4]=e_3
\end{array}\right.&
L(\alpha,0): \left\{\begin{array}{l} [e_{i},e_1]=e_{i+1},\quad 1\leq i\leq 2\\{}
[e_1,e_4]=\alpha e_2+f_2\\{}
 [e_2,e_4]=\alpha e_3
 \end{array}\right.
\end{array}$$

If we considered algebra
$$L(\alpha,0): \left\{\begin{array}{ll}
[e_{i},e_1]=e_{i+1}&1\leq i\leq 2\\{} [e_1,e_4]=\alpha e_2+f_2\\{}
 [e_2,e_4]=\alpha e_3
 \end{array}\right.
$$

We have $\alpha=1$ or $0$. And we obtain the two first algebras of proposition. By
standard way it is not difficult to check that these algebras are not isomorphic.

Consider algebra $L(\alpha,1)$, we make the general change of basis
$$e'_1=a_1e_1+a_2e_4, \ \ e'_4=b_1e_1+b_2e_4$$ where $a_1b_2-a_2b_1\neq 0.$

In other hand $[e'_4,e'_1]=0$ and we have
$$b_1a_1+b_1a_2\alpha=0$$
$$b_1a_2=0$$
it implies that $b_1=0$. Finally we obtain
$$\alpha'=\frac{b_2[a_1\alpha+a_2(\alpha^2+1)]}
{[(a_1+a_2\alpha )^2+a_2^2]}$$

Comparing the coefficients at the basic element we obtain
restriction
$$b_2^2=\frac{[(a_1+a_2\alpha )^2+a_2^2]^2} {a_1^2}$$

It is not difficult to check that the nullity of the following
expression is invariant because:
$$1+\alpha'^2=\frac{(1+\alpha ^2)((a_1+a_2 \alpha)^2+a_2^2)}
{a_1^2}=$$

\textbf{Case 1.} $\alpha^2+1 \neq 0$ then putting
$a_2=-\frac{a_1\alpha} {1+\alpha^2}$ implies $\alpha'=0$. Thus, in
this case we obtain $\mu_{4}$.

\

\textbf{Case 2.} $\alpha^2+1 = 0$ (i.e $\alpha=\pm i$) then we
have that $b_2=\pm \frac{(a_1+a_2\alpha^2)+a_2^2}{a_1}$ and
$\alpha'=\pm \alpha$ we obtain $\alpha'=i$. Thus, in this case we
obtain $\mu_3$. $\hfill\Box$

\

\begin{thm} \label{(I,1,2)}
Let $\ll$ be an $n$-dimensional ($n\geq 6$) graded $2$-filiform non split Leibniz algebra of type
$\mu_{(I, 1, 2)}$. Then $\ll$ is isomorphic to the one of the following pairwise non
isomorphic algebras:
$$\begin{array}{ll}
 \left\{\begin{array}{ll}
[e_i, e_1]=e_{i+1}, &  1\leq i\leq n-3 \\[1mm]
[e_1, e_{n-1}]=e_2+e_n & \\[1mm]
[e_i, e_{n-1}]=e_{i+1}, &  2\leq i\leq n-3
\end{array}\right.&
 \left\{\begin{array}{ll}
[e_i, e_1]=e_{i+1}, &  1\leq i\leq n-3 \\[1mm]
[e_1, e_{n-1}]=e_n &
\end{array}\right.
\end{array}$$
\end{thm}

{\bf Proof:}

According to the theorem conditions we have the following
multiplications in an adapted basis $\{e_1, e_2,\dots, e_n\}$:

$$\left\{\begin{array}{ll}
[e_i, e_1]=e_{i+1}, &  1\leq i\leq n-3 \\[1mm]
[e_1, e_{n-1}]={\alpha_1}e_2+ {\gamma_1}e_n& \\[1mm]
[e_i, e_{n-1}]={\alpha_i}e_{i+1}, &  2\leq i\leq n-3  \\[1mm]
[e_{n-1}, e_{n-1}]={\alpha_{n-1}}e_2+{\gamma_{n-1}}e_n &\\[1mm]
[e_n, e_{n-1}]={\alpha_n}e_3 &\\[1mm]
[e_i, e_n]={\beta}_ie_{i+2}, & 1\leq i\leq n-4 \\[1mm]
[e_{n-1}, e_n]={\beta_{n-1}}e_3& \\[1mm]
[e_n, e_n]={\beta_n}e_4&
\end{array}\right.$$
where either $\gamma_1 \neq 0$ or $\gamma_{n-1} \neq 0$.

Using Leibniz identity it is not difficult to obtain the following restrictions:
$$\left\{\begin{array}{ll}
\alpha_i=\alpha, & 1\leq i\leq n-3\\
{\beta_i}{\gamma_1}=0, & 1\leq i\leq n-4\\
{\beta_i}\gamma_{n-1}=0,& 1\leq i\leq n-4\\
{\gamma_1}{\beta_{n-1}}+\alpha_{n-1}=0\\
\alpha_{n-1}=\alpha_n=0&\\
\beta_i\gamma_j=0, & i\in\{n-1,n\}, \ j\in\{1,n-1\}
\end{array}\right.$$

Since either $\gamma_1 \neq 0$ or $\gamma_{n-1} \neq 0$, we have
that  $\beta_i=0$ for $1 \leq i \leq n-4$ and $\beta_{n-1}=\beta_n=0$. Thus, the
multiplications have the following form:
$$\left\{\begin{array}{ll}
[e_i, e_1]=e_{i+1}, &  1\leq i\leq n-3 \\[1mm]
[e_1, e_{n-1}]={\alpha}e_2+ {\gamma_1}e_n& \\[1mm]
[e_i, e_{n-1}]={\alpha}e_{i+1}, &  2\leq i\leq n-3  \\[1mm]
[e_{n-1}, e_{n-1}]={\gamma_{n-1}}e_n&
\end{array}\right.$$

It is possible to suppose that
$$R_{e_1+Ae_{n-1}}=\left(\begin{array}{ll}
(1+A\alpha)I_{n-3}& \begin{array}{lll}
0&0&0\\
\vdots&\vdots&\vdots\\
0&0&0
\end{array}  \\
\begin{array}{lll}
0 &0&\cdots\\
0&0&\cdots\\
A\gamma_1 &0& \cdots
\end{array}&
\begin{array}{lll}
0&0&0\\
0&0&0\\
0&0&A\gamma_{n-1}
\end{array}
\end{array}\right)$$
where $I_{n-3}$ is the unit matrix of size $n-3$ and $1+A\alpha\neq 0$.

As $rang(R_{e_1+Ae_{n-1}})\leq n-3$ (otherwise the characteristic
sequence for element $e_1+A e_{n-1}$ would be greater than $(n-p, 1, \dots,
1)$), then $(1+A\alpha)^{n-3} A\gamma_{n-1}=0$, hence
$\gamma_{n-1}=0$ and $\gamma_1\neq 0$. By an elementary change of
basis, it is possible to suppose that $\gamma_1=1$.

By a general change of basis the expression for the new generators is
$$e'_1=\displaystyle\sum_{i=1}^{n-1} A_i e_i,\qquad \qquad
e'_{n-1}=\displaystyle\sum_{i=1}^{n-1} B_i e_i$$
obtaining $\alpha'=\displaystyle\frac{B_{n-1}\alpha}{A_1+A_{n-1}\alpha}$.

It is easy to see that if $\alpha\neq 0$ we have the first algebra of the theorem
and if $\alpha=0$ the second algebra.

$\hfill\Box$

\

\

Consider now graded $2$-filiform Leibniz algebras of type II.

\

Let $\ll$ be a graded $n$-dimensional $p$-filiform Leibniz
algebra. Then there exists a basis $\{e_1, e_2,\dots, e_{n-p},f_1,f_2,\dots,f_p\}$ of
$\ll$ such that multiplications on the right side on element $e_1$
will have the form:
$$\left\{\begin{array}{ll}
[e_1,e_1]=0&\\[1mm]
[e_i,e_1]=e_{i+1}, &  2\leq i\leq n-p-1 \\[1mm]
[f_j,e_1]=0, &  1\leq j\leq p
\end{array}\right.$$

From these multiplications we have: $$<e_1>\subseteq \ll_1,\qquad
<e_{i+1}>\subseteq \ll_i,\quad 2 \leq i \leq n-2$$

But again we do not know about the position of elements $\{e_2, f_{2},
f_{3},\dots,f_p\}$ in natural gradation.

Let denote by $r_1, r_2,\dots,r_p$ $(r_1\leq r_2 \leq \cdots \leq r_p)$ the places of elements
$e_2,$ $ f_{2},$ $ f_{3},$ $\dots,$ $ f_p$
correspondingly, that is, $e_2\in \ll_{r_1}$, $f_{i}\in \ll_{r_i}$,  $2\leq i \leq p$.

Let $\ll$ be a graded $2$-filiform Leibniz algebra. Since $r_1=1$, further we shall denote
$r_2$ by $r$.

For the $2$-filiform Leibniz algebras of type II the following lemma is hold.

\begin{lem}
Let $\ll$ be an $n$-dimensional
$2$-filiform Leibniz algebra. Then the
following conditions are hold:
\begin{itemize}
\item[$a)$] $\ll$ has nilindex $n-1$;\\[2mm]
\item[$b)$] or $dim(\ll^i)=n-1-i$,\quad $2 \leq i \leq n-2$

 \noindent or $dim(\ll^i)=\left\{\begin{array}{ll}
n-i,&2\leq i\leq r\\
n-1-i~,&r+1\leq i\leq n-2
\end{array}\right.$
for some  $r$,  $2\leq r\leq n-2$
\end{itemize}
\end{lem}

{\bf Proof:}

$a)$ Let $x\in \ll-[\ll,\ll]$ such that $C(x)=(n-2, 1,1)$.
Hence, $R_x^{n-2}=0$ and $R_x^{n-3}\neq 0$ and, consequently,
there exists element $y\in\ll$, such that $R_x^{n-3}(y)\neq 0$.
Therefore $\ll^{n-2}\neq 0$ and $\ll^{n-1}=0$ (when $\ll^{n-1}\neq 0$,
then by [\cite{Omirov1},lemma 1, lemma 4] the algebra $\ll$ would be either nul-filiform or filiform).

\

$b)$ Let $e_1\in \ll-[\ll,\ll]$ be a maximal characteristic vector of $ll$, where $\ll$ is of type I.
Then for $r=1$, that is, $dim(\ll/\ll^2)=3$ we have that $dim(\ll^i)=n-1-i$, $2 \leq i\leq n-2$.
For $r=2$, that is, $dim(\ll/\ll^2)=2$ we get:
$$dim(\ll^i)=\left\{\begin{array}{ll}
n-2,&i=2\\
n-1-i,&3\leq i\leq n-2
\end{array}\right.$$

Let algebra $\ll$ has the type II. For $r_2=1$ we obtain that $dim(\ll/\ll^2)=3$, that is,
$dim(\ll^i)=n-1-i$, $2 \leq i \leq n-2$. For $r_2\in \{2, 3, \dots, n-2\}$ we get:
$dim(\ll/\ll^2)=2$, that is, $dim(\ll^i)=\left\{\begin{array}{ll}
n-i,&2\leq i\leq r\\
n-1-i,&r+1\leq i\leq n-2
\end{array}\right.$
$\hfill\Box$

\begin{lem}\label{typeII,r>2}
Let $\ll$ be a complex n-dimensional $(n \geq 5)$ graded
$2$-filiform Leibniz algebra of type II and $r > 2$. Then $\ll$ is a Lie
algebra.
\end{lem}

{\bf Proof:}

Let $(1)$ be the family of laws of $\ll$:
$$(1)\left\{\begin{array}{ll}
[e_i, e_1]=e_{i+1}, &  2\leq i\leq n-2 \\[1mm]
[e_1, e_i]={\alpha_{1,i}}e_{i+1}, & 2\leq i\leq n-2, \ i\neq r \\[1mm]
[e_1, e_r]={\alpha_{1,r}}e_{r+1}+{\gamma_1}e_n & \\[1mm]
[e_1, e_n]={\alpha}_{1,n}e_{r+2} &\\[1mm]
[e_i, e_j]={\alpha_{i,j}}e_{i+j-1}, & 2\leq i,j\leq n-2, \ i+j\leq n, \ i+j\neq r+2 \\[1mm]
[e_i, e_{r+2-i}]={\alpha}_{i,r+2-i}e_{r+1}+{\gamma_i}e_n, & 2\leq i\leq r  \\[1mm]
[e_n, e_i]={\alpha_{n,i}}e_{i+r}, & 2\leq i\leq n-r-1 \\[1mm]
[e_i, e_n]={\alpha_{i,n}}e_{i+r}, & 2\leq i\leq n-r-1 \\[1mm]
[e_n, e_n]={\alpha_{n,n}}e_{2r+1}, & r\leq \frac{n-2}{2}
\end{array}\right.$$
where omitted products are zero and $(\gamma_1,\gamma_2,\dots,\gamma_r)\neq (0,0,\dots,0)$.

Using Leibniz identity we get the following restrictions:
$$\left\{\begin{array}{ll}
\alpha_{1,i}=\alpha, & 2\leq i\leq n-2\\
\gamma_1=0&\\
{\alpha_1}(\alpha_1+1)=0&\\
\alpha_{1,n}=0,& r\leq n-4\\
\alpha_{n,n}=0,
& r\leq \displaystyle\frac{n-3}{2}
\end{array}\right.$$

It is necessary to consider separately the cases $r=n-3$, $r=n-2$ and $r=\frac{n-2}{2}$ ($n$ even).

{\bf Case 1}. $\alpha=0$. Then $(1)$ will have the following form:
$$\left\{\begin{array}{ll}
[e_i, e_1]=e_{i+1}, &  2\leq i\leq n-2 \\[1mm]
[e_i, e_j]={\alpha_{i,j}}e_{i+j-1}, & 2\leq i,j\leq n-2, \ i+j\leq n, \ i+j\neq r+2 \\[1mm]
[e_i, e_{r+2-i}]={\alpha}_{i,r+2-i}e_{r+1}+{\gamma_i}e_n, & 2\leq i\leq r  \\[1mm]
[e_n, e_i]={\alpha_{n,i}}e_{i+r}, & 2\leq i\leq n-r-1 \\[1mm]
[e_i, e_n]={\alpha_{i,n}}e_{i+r}, & 2\leq i\leq n-r-1
\end{array}\right.$$

Using Leibniz identity for elements $\{e_i, e_{r+1-i}, e_1\}$ for
$2\leq i \leq r$ and $\{e_i,e_1,e_{r+1-i}\}$ for $2\leq i\leq r$, that is, $$[e_i, [e_{r+1-i}, e_1]]=[[e_i,
e_{r+1-i}], e_1] - [[e_i, e_1],e_{r+1-i}]$$
$$[e_i, [e_1,e_{r+1-i}]]=[[e_i,e_1],e_{r+1-i}] - [[e_i,e_{r+1-i}], e_1]$$
we obtain that
$\gamma_i=0$ for $2 \leq i \leq r$. Hence $e_n \notin \ll^2$
and $r=1$ we have the contradiction to the condition of the lemma.

{\bf Case 2}. $\alpha=-1$. Then the multiplications $(1)$ will have the form:
$$\left\{\begin{array}{ll} [e_i,e_1]=e_{i+1},&  2\leq i\leq n-2
\\{} [e_1, e_i]=-e_{i+1},&2\leq i\leq n-2 \\{} [e_i,
e_j]={\alpha_{i,j}}e_{i+j-1}, & 2\leq i,j\leq n-2, \quad i+j\leq n
\quad i+j\neq r+2 \\{} [e_i,
e_{r+2-i}]={\alpha_{i,r+2-i}}e_{r+1}+{\gamma_i}e_n, & 2\leq i\leq
r \\{} [e_n, e_i]={\alpha_{n,i}}e_{i+r}, & 2\leq i\leq n-r-1 \\{}
[e_i, e_n]={\alpha_{i,n}}e_{i+r}, & 2\leq i\leq n-r-1
\end{array}\right.$$
where $(\gamma_2,\dots,\gamma_r)\neq (0,\dots,0)$.

Using Leibniz identity it is not difficult to get that
$$\begin{array}{ll}
\alpha_{n,i}=\alpha_n, & 2\leq i\leq n-r-1\\
\alpha_{i,n}=\overline{\alpha_n}, & 2\leq i\leq
n-r-1\\
\alpha_n=-\overline{\alpha_n}&
\end{array}$$

From equality $[e_1, [e_i, e_i]]=0$ for $2\leq i\leq
\frac{n-3}{2}$, we have $\alpha_{i,i}=0$ for $2\leq i\leq
\frac{n-1}{2}$.

When $i=\frac{n}{2}$ (when $n$ is even), we consider the
following equalities:
\begin{enumerate}
\item[$ $] $[e_{\frac{n}{2}}[e_{\frac{n}{2}-1},e_1]]=[[e_{\frac{n}{2}},e_{\frac{n}{2}-1}],e_1]-
[[e_{\frac{n}{2}},e_1],e_{\frac{n}{2}-1}]$\hfill $(2)$

\item[$ $] $[e_{\frac{n}{2}-1}[e_{\frac{n}{2}-1},e_1]]=[[e_{\frac{n}{2}-1},e_{\frac{n}{2}-1}],e_1]-
[[e_{\frac{n}{2}-1},e_1],e_{\frac{n}{2}-1}]$\hfill $(3)$

\item[$ $] $[e_1,[e_{\frac{n}{2}-1},e_{\frac{n}{2}}]]=[[e_1,e_{\frac{n}{2}-1}],e_{\frac{n}{2}}]-
[[e_1,e_{\frac{n}{2}}],e_{\frac{n}{2}-1}]$\hfill $(4)$
\end{enumerate}

From equalities $(2) $ up to $(4)$ we obtain the restrictions:
$$(5)\left\{\begin{array}{l}
\alpha_{\frac{n}{2},\frac{n}{2}}=\alpha_{\frac{n}{2},\frac{n}{2}-1}-\alpha_{\frac{n}{2}+1,\frac{n}{2}-1}\\[2mm]
\alpha_{\frac{n}{2},\frac{n}{2}-1}=-\alpha_{\frac{n}{2}-1,\frac{n}{2}}\\[2mm]
\alpha_{\frac{n}{2},\frac{n}{2}}=\alpha_{\frac{n}{2}+1,\frac{n}{2}-1}+\alpha_{\frac{n}{2}-1,\frac{n}{2}}
\end{array}\right.$$

From $(5)$ we have that $\alpha_{\frac{n}{2},\frac{n}{2}}=0$.

Thus, we prove that $[e_i,e_i]=0$ for $1\leq i \leq n$.

From the following chain of equalities:
$$\begin{array}{lll}
[e_i,e_j]&=[e_i,[e_{j-1},e_1]]=[[e_i,e_{j-1}],e_1]-[[e_i,e_1],e_{j-1}]=\\
&=-[e_1,[e_i,e_{j-1}]]+[[e_1,e_i],e_{j-1}]=\\
&=-([[e_1,e_i],e_{j-1}]-[[e_1,e_{j-1}],e_i])+[[e_1,e_i],e_{j-1}]=[[e_1,e_{j-1}],e_i]=\\
&=-[e_j,e_i]
\end{array}$$
we obtain that $[e_i,e_j]= - [e_j, e_i]$ for $1\leq i < j \leq n$, that is, is a Lie algebra.

The cases $r=n-3$, $r=n-2$ and $r=\frac{n-2}{2}$  (when $n$ is even) are proved analogously.
$\hfill\Box$

\

\

Next, we will see some examples of graded filiform Leibniz algebras of type II.

\

{\bf Example 2}. Let $\ll$ be a graded filiform Leibniz algebra
of type II. Then $\ll\oplus\CC$ is graded $2$-filiform Leibniz algebra of type
II.

\

And now, we prove some lemmas for a graded non split and non Lie $2$-filiform
Leibniz algebra of type II.

\begin{lem}\label{typeII,r=1}
There exits no a graded non split and non Lie $2$-filiform Leibniz algebra of
type II and $r=1$.
\end{lem}

{\bf Proof:}

Let $\ll$ be a Leibniz algebra which satisfies
the condition of the lemma. Then the table of multiplication is

$$\left\{\begin{array}{ll}
[e_i, e_1]=e_{i+1}, &  2\leq i\leq n-2 \\[1mm]
[e_1, e_i]={\alpha_{1,i}}e_{i+1}, & 2\leq i\leq n-2 \\[1mm]
[e_1, e_n]={\alpha_{1,n}}e_3 & \\[1mm]
[e_i, e_j]={\alpha_{i,j}}e_{i+j-1}, &2\leq i,j\leq n-2, \quad i+j\leq n \\[1mm]
[e_n, e_i]={\alpha_{n,i}}e_{i+1}, & 2\leq i\leq n-2 \\[1mm]
[e_i, e_n]={\alpha_{i,n}}e_{i+1}, & 2\leq i\leq n-2 \\[1mm]
[e_n, e_n]={\alpha_{n,n}}e_3&
\end{array}\right.$$

From Leibniz identity we have the following restrictions:
$$\begin{array}{ll}
\alpha_{1,i}=\alpha, & 2\leq i\leq n-2\\
{\alpha_1}({\alpha_1}+1)=0&\\
\alpha_{1,n}=\alpha_{n,n}=0
\end{array}$$

{\bf Case 1}. $\alpha=0$. Using Leibniz identity we get
$$\begin{array}{ll}
\alpha_{i,j}=\alpha_j, & 2\leq i\leq n-2\\
\alpha_j=0,&3\leq j\leq n-2\\
{\alpha_{i,n}}=\alpha_n&2\leq i\leq n-2\\
\alpha_{n,i}=0&2\leq i\leq n-2
\end{array}$$
and taking a change of basis: $e'_2=e_2-\alpha_2 e_1,$\
$e'_n=e_n-\alpha_n e_1,$\ $e'_i=e_i$ for $i\neq 2, n$, we have
that  $\alpha_2=\alpha_n=0$, that is, $\ll$ is split.

{\bf Case 2}. $\alpha =-1$. Analogous to lemma \ref{typeII,r>2}, we get a Lie algebra.

$\hfill\Box$

\

\

\begin{lem}\label{typeII,r=2}
There exits no a graded non split and non Lie $2$-filiform Leibniz algebra of type II and $r=2$.
\end{lem}

{\bf Proof:}

Let $\ll$ be a Leibniz algebra satisfying the conditions of
the lemma. Then, there exists an adapted basis $\{e_1, e_2,\dots, e_n\}$ of $\ll$
such that the multiplications will be the following:
$$\left\{\begin{array}{ll}
[e_i, e_1]=e_{i+1},& 2\leq i\leq n-2 \\[1mm]
[e_1, e_2]={\alpha_{1,2}}e_3+{\gamma_1}e_n \\[1mm]
[e_1, e_i]={\alpha_{1,i}}e_{i+1}, & 3\leq i\leq n-2 \\[1mm]
[e_1, e_n]={\alpha}_{1,n}e_4 \\[1mm]
[e_2, e_2]={\alpha_{2,2}}e_3+{\gamma_2}e_n \\[1mm]
[e_i, e_j]={\alpha_{i,j}}e_{i+j-1}, & 2\leq i,j\leq n-2, \ i+j\leq n, \ (i,j)\neq (2,2) \\[1mm]
[e_n, e_i]={\alpha_{n,i}}e_{i+2}, &  2\leq i\leq n-3 \\[1mm]
[e_i, e_n]={\alpha_{i,n}}e_{i+2}, &  2\leq i\leq n-3 \\[1mm]
[e_n, e_n]={\alpha_{n,n}}e_5, &
\end{array}\right.$$
where either $\gamma_1\neq 0$ or $\gamma_2\neq 0$.

As in the above two lemmas we obtain:
$$\left\{\begin{array}{ll}
\alpha_{1,i}=\alpha,& 2\leq i\leq n-2\\
{\alpha}(1+\alpha)=0&\\
\alpha_{1,n}=\alpha_{n,n}=0&
\end{array}\right.$$

Let us consider two possible cases for parameter $\alpha$.

{\bf Case 1}.  $\alpha=0$. Then, the multiplications in $\ll$ have the form:
$$\left\{\begin{array}{ll}
[e_i, e_1]=e_{i+1},& 2\leq i\leq n-2 \\[1mm]
[e_1, e_2]={\gamma_1}e_n \\[1mm]
[e_2, e_2]={\alpha_{2,2}}e_3+{\gamma_2}e_n \\[1mm]
[e_i, e_j]={\alpha_{i,j}}e_{i+j-1}, & 2\leq i,j\leq n-2, \ i+j\leq n, \ (i,j)\neq (2,2) \\[1mm]
[e_n, e_i]={\alpha_{n,i}}e_{i+2}, &  2\leq i\leq n-3 \\[1mm]
[e_i, e_n]={\alpha_{i,n}}e_{i+2}, &  2\leq i\leq n-3
\end{array}\right.$$
where either $\gamma_1\neq 0$ or $\gamma_2\neq 0$.

Using Leibniz identity leads us to the following restrictions:
$$\begin{array}{ll}
\alpha_{i,j}=\alpha_j, & 2\leq i,j\leq n-2\\
\alpha_j=0, & \ 3\leq j\leq n-2\\
\alpha_{i,n}=\alpha_n, & 2\leq i\leq n-3\\
\alpha_{n,i}=0, & 2\leq i\leq n-3\\
\alpha_n \gamma_2=0&\\
\alpha_n\gamma_1=0&
\end{array}$$

Either $\gamma_1\neq 0$ or $\gamma_2\neq 0$ (otherwise algebra $\ll$ is split), then $\alpha_n=0$.

The change of basis given by $e'_2=e_2-{\alpha_2}e_1,$ $e'_i=e_i,
\ i\neq 2,$ allows to suppose $\alpha_2=0$.

Consider the operator of right multiplication $R_{e_1+Ae_{n-1}}$,
where $0\neq A\in \CC$ such that $e_1+A e_{n-1}\neq 0$. $\gamma_2=0$ and hence $\gamma_1\neq 0$
may be proved in much the same way as the proof of theorem \ref{(I,1,2)}.
Without loss of generality we can assume that $\gamma=1$.

Thus, we have the following multiplications in algebra $\ll$:

$$\left\{\begin{array}{ll}
[e_i, e_1]=e_{i+1},& 2\leq i\leq n-2 \\[1mm]
[e_1, e_2]=e_n &
\end{array}\right.$$

Taking the change of basis in form:
$$\begin{array}{l}
e'_1=e_1+e_2,\quad e'_2=e_3+e_n,\\
e'_i=e_{i+1}, \quad 3\leq i\leq n-2,\\
e'_{n-1}=e_1,\quad e'_n=e_n
\end{array}$$
we obtain the algebra of type I.

{\bf Case 2}. $\alpha=-1$. As in above cases, we get a Lie algebra.

$\hfill\Box$

\

\

From lemmas \ref{typeII,r>2}, \ref{typeII,r=1} and \ref{typeII,r=2} we can conclude that there exist no
graded non split and non Lie $2$-filiform  Leibniz algebras of type II.

\

Thus, according to theorem \ref{(I,1,2)} we have the classification of non split and non
Lie $2$-filiform Leibniz algebras. Summing the classification of non split
graded $2$-filiform Lie algebras \cite{AJM} and the result of theorem \ref{(I,1,2)}
we have completed the classification of graded non split $2$-filiform
Leibniz algebras.

\mbox{}

\end{document}